\documentclass[preprint,12pt]{elsarticle}

\usepackage{natbib}
\usepackage{amsmath,amsfonts,amssymb}

\usepackage{graphicx}
\usepackage{caption}
\usepackage{subcaption}

\newtheorem{proposition}{Proposition}

\newtheorem{theorem}{Theorem}
\newtheorem{lemma}{Lemma}
\newtheorem{example}{Example}
\newtheorem{remark}{Remark}
\newtheorem{const}{Construction}

\newproof{pf}{Proof}
\newproof{pot}{Proof of Theorem \ref{Th1}}

\def\primal{\mathrm{primal}}
\newcommand{\R}{\mathbb{R}}
\newcommand{\Q}{\mathbb{Q}}
\newcommand{\D}{\mathbb{D}}
\newcommand{\Dh}{\mathbb{DH}}
\newcommand{\qi}{\mathbf{i}}
\newcommand{\qj}{\mathbf{j}}
\newcommand{\qk}{\mathbf{k}}

\newcommand{\eps}{\epsilon}

\journal{Mechanism and Machine Theory}

\begin{document}

\begin{frontmatter}

\title{Classification of Angle-Symmetric 6R Linkages}

\author{Zijia Li and Josef Schicho\corref{cor1}}
\ead{\{zijia.li,\ josef.schicho\}@oeaw.ac.at}
\cortext[cor1]{Corresponding author.}
\address{Johann Radon Institute for
  Computational and Applied Mathematics Austrian Academy of Sciences
  (RICAM), Altenbergerstrasse 69, 4040 Linz, Austria}

\begin{abstract}
In this paper, we consider a  special kind of overconstrained 6R closed linkages which
  we call angle-symmetric 6R linkages. These are linkages with the property that the rotation angles are
 equal for each of the three pairs of opposite joints. 
 We give a classification of these
 linkages. It turns that there are three types. First, we have the linkages with line symmetry. The second type is new.
 The third type is related to cubic motion polynomials.
\end{abstract}


\end{frontmatter}

\section{Introduction}
\label{intr}
Movable closed 6R linkages have been considered by many authors (see~\cite{Sarrus, Goldberg, Waldron, B80, Wohlhart91, Diet95}).
 In this paper, we give the complete classification of a certain class of such linkages,
 which we call angle-symmetric. This means that the rotation angles at the three
 pairs of opposite joints are equal for all possible configurations, or at least
 for infinitely many configurations (it could be that a certain linkage has two
 components, where only one of them is angle-symmetric). It is well-known that
 the line symmetric linkage of Bricard \cite{B80} is angle-symmetric.  A second
 family is new; it can be characterized by the presence of three pairs of parallel
 rotation axes. This fills a gap in~\cite[Section 3.8]{B03}.
 A third family 
 was discovered in \cite{HSScubep, part11} using factorizations of cubic motion polynomials.

Our main tool is the $\lambda$-matrix of a linkage, to be defined in section {\ref{lmatrix}},
 and its rank $r$. Intuitively speaking, the configuration set can be described as
 the vanishing set of $r$ equations in three variables, namely the cotangents of the half
 of the rotation angles. We will show that $r$ is either 2, 3, or 4. If $r=2$,
 then the linkage is line symmetric. If $r=3$, then we get the new linkage
 with three pairs of parallel axes. If $r=4$, then we obtain the linkage described
 in \cite{HSScubep, part11} using motion polynomials.

We use Study's description of Euclidean displacements by the algebra $\mathbb{DH}$ 
 of dual quaternions 
 (see \cite{HSScubep, part11}).

\paragraph*{Structure of the paper} The remaining part of the paper is set up as follows. In Section 2,
 we give the definition of the $\lambda$-matrix. 
 We also show that the rank of this matrix is 2, 3, or 4. 
 Section 3 contains the main result and examples

\section{The $\lambda$-matrix}\label{lmatrix}

In this section we define, for a given linkage, a matrix ${\bf M_\lambda}$ whose rows
 are related to an algebraic system defining the configuration space. In the next
 section, we will see that the rank of this matrix is the basic criterion for
 classifying angle-symmetric linkages.

 The set of all possible motions of a closed 6R linkage is determined by the
 position of the six rotation axes in some fixed initial configuration. 
 (The choice of the initial configuration among all possible configurations
 is arbitrary. In some later steps in the classification, we will occasionally change the
 initial configuration.)

The algebra $\mathbb{DH}$ of dual quaternions is the 8-dimensional real vector space
 generated by $1,\eps,\qi,\qj,\qk,\eps\qi,\eps\qj,\eps\qk$ (see \cite{HSScubep, part11}).
 Following \cite{HSScubep, part11}, we can represent a rotation by a dual quaternion of the
 form $\left(\cot\left(\frac{\phi}{2}\right)-h\right)$, where $\phi$ is the rotation
 angle and $h$ is a dual quaternion such that $h^2=-1$ depending only on the rotation axis.
 We use projective representations, which means that two dual quaternions represent
 the same Euclidean displacement if only if one is a real scalar multiple of the other.

Let $L$ be a 6R linkage given by 6 lines, represented by dual quaternions 
 $h_1,\dots,h_6$ such that $h_i^2=-1$ for $i=1,\dots,6$. A configuration (see \cite{HSScubep, part11}) is a 6-tuple
 $(t_1,\dots,t_6)$, such that the closure condition
 $$(t_1-h_1)(t_2-h_2)(t_3-h_3)(t_4-h_4)(t_5-h_5)(t_6-h_6)\in \mathbb{R}\backslash \{0\}$$
 holds. The configuration parameters $t_i$ -- the cotangents of the rotation angles --
 may be real numbers or $\infty$, and in the second case we evaluate the expression
 $(t_i-h_i)$ to $1$, the rotation with angle 0.
 The set of all configurations of $L$ is denoted by $K_L$.
 
 There is a subset of $K_L$, denoted by $K_{sym}$, defined by the additional restrictions 
 $t_1=t_4,t_2=t_5,t_3=t_6$. We assume that $K_{sym}$ is a one-dimensional set, i.e. the linkage 
 has an angle-symmetric motion. Mostly, we will assume, slightly stronger, that there exists an irreducible 
 one-dimensional set for which none of the $t_i$ is fixed. Such a component is called a non-degenerate component.
 We also exclude the case  $\dim_{\mathbb{C}}K_{sym}\geq 2$. 
 Linkages with mobility $\geq 2$ do exist, but they are well understood.

 The closure condition is equivalent to
\[(t_1-h_1)(t_2-h_2)(t_3-h_3)=\lambda (t_3+h_6)(t_2+h_5)(t_1+h_4),\]
 where $\lambda$ is a nonzero real value depending on $t_1,t_2,t_3$.
 By taking norm on both sides, we get $\lambda^2=1$, i.e. $\lambda=\pm 1$. 
 By multiplying both sides with $(t_1+h_1)$ from the left and with $(t_1-h_4)$
 from the right, and afterwords dividing by $(t_1^2+1)$, we obtain the equation
 \[(t_2-h_2)(t_3-h_3)(t_1-h_4)= \lambda (t_1+h_1)(t_3+h_6)(t_2+h_5).\]
 Similarly, we obtain
\[(t_3-h_3)(t_1-h_4)(t_2-h_5)= \lambda (t_2+h_2)(t_1+h_1)(t_3+h_6),\]
\[(t_1-h_4)(t_2-h_5)(t_3-h_6)= \lambda (t_3+h_3)(t_2+h_2)(t_1+h_1),\]
\[(t_2-h_5)(t_3-h_6)(t_1-h_1)= \lambda (t_1+h_4)(t_3+h_3)(t_2+h_2),\]
\[(t_3-h_6)(t_1-h_1)(t_2-h_2)= \lambda (t_2+h_5)(t_1+h_4)(t_3+h_3).\]

We may divide $K_{sym}$ into two disjoint subsets $K_{sym}^+$ and $K_{sym}^-$, according to whether $\lambda$
 is equal to $+1$ or $-1$ in the equations above. Any irreducible component of $K_{sym}$ is either contained in $K_{sym}^+$
 or in $K_{sym}^-$. Note that $\infty^3$ is an element of $K_{sym}^+$.

\begin{remark}\label{rm2}
 When we want to study some component $K_0\subset K_{sym}^-$, we may proceed in the following
 way: we take a configuration $\tau \in K_0$, which defines  a set of  rotations  around the joint axes.
 Then we apply these rotations, obtaining new positions for the 6 lines. In the transformed linkage, the component corresponding
 to $K_0$ contains $\infty^3$. So we will always assume that $ \lambda=1$.
\end{remark}

When $\lambda=1$, after moving  the right parts of the above equations to the left, we get an equation
\begin{equation*}\label{I02}
\bf{ M^{\dagger} } \bf{X}=\bf{0},             
\end{equation*}
where  ${\bf X}=\left[
      {t_1 t_2},     
      {t_1 t_3},     
      {t_2 t_3},     
          {t_3},     
          {t_2},     
          {t_1},     
            {1}     
 \right]^{T}$. 
 If we denote $h_6+h_3, h_5+h_2, h_4+h_1$ by $g_3,g_2,g_1$ respectively, then the coefficient 
 matrix ${\bf M^{\dagger} }$ is 
 \begin{equation*}\label{I03}
\left[
 \begin{array}{ccccccc}
 g_3, g_2, g_1,& h_5h_4-h_1h_2,& h_6h_4-h_1h_3,& h_6h_5-h_2h_3,& h_6h_5h_4+h_1h_2h_3 \\
 g_3, g_2, g_1,& h_1h_5-h_2h_4,& h_1h_6-h_3h_4,& h_6h_5-h_2h_3,& h_1h_6h_5+h_2h_3h_4 \\
 g_3, g_2, g_1,& h_2h_1-h_4h_5,& h_1h_6-h_3h_4,& h_2h_6-h_3h_5,& h_2h_1h_6+h_3h_4h_5 \\
 g_3, g_2, g_1,& h_2h_1-h_4h_5,& h_3h_1-h_4h_6,& h_3h_2-h_5h_6,& h_3h_2h_1+h_4h_5h_6 \\
 g_3, g_2, g_1,& h_4h_2-h_5h_1,& h_4h_3-h_6h_1,& h_3h_2-h_5h_6,& h_4h_3h_2+h_5h_6h_1 \\
 g_3, g_2, g_1,& h_5h_4-h_1h_2,& h_4h_3-h_6h_1,& h_5h_3-h_6h_2,& h_5h_4h_3+h_6h_1h_2 \\
       \end{array}\right].          
\end{equation*}
Note that $\bf{ M^{\dagger} }$ is a $6\times7$ matrix with entries in dual quaternions. We also consider 
 $\bf{ M^{\dagger} }$ to be a $48\times7$ matrix with real entries.
 It can be decomposed into submatrices $M^{\dagger}_1,\dotsm,M^{\dagger}_6$, 
 where $M^{\dagger}_i$ is the real $8\times7$ matrix -- or the row vector with 7 dual quaternion entries --
 corresponding to the $i-th$ equivalent formulation of the closure condition above, for
 $i=1,\ldots,6$.

 Our classification is based on the following theorem which gives the bounds for the rank 
 of $\bf{ M^{\dagger} }$.
 \begin{theorem}\label{Th1}
 Assume that ${K_{sym}}$ contains a non-degenerate component of dimension $1$. 
 Then $r:=\operatorname{rank}(\bf{ M^{\dagger} })$ $ \in \{ 2,3,4\}$.
 \end{theorem}

Before we prove Theorem \ref{Th1}, we give a lemma.

\begin{lemma}\label{lm1}
 Assume that ${K_{sym}}$ contains a non-degenerate component $K_0$ of dimension $1$ such that
 $\infty^3\in K_0$, and $r\geq 4$.  
 Then there exists a polynomial of the form  
 $$b t_1 +c t_2 +d,$$
where $b,c,d \in \mathbb{R}$ and $bc \neq 0$,
 which vanishes on  ${K_{sym}}$,
 maybe after some permutation of the variables $t_1,t_2,t_3$.
 Moreover, we can define a matrix $\bf{ N^{\dagger} }$ of rank $\geq r-2$ such that the projection of ${K_{sym}}$
 to $(t_1,t_3)$ is defined by
 \begin{equation}
  \bf{ N^{\dagger} } \bf{X'}=0,
 \end{equation}
 where ${\bf{X'}}=[t_1^2,t_1 t_3,t_1,t_3,1]^{T}$.
 \end{lemma}
\begin{pf}
  As $r\geq 4$, we have at least four independent equations in three variables 
 $(t_1,t_2,t_3)$ of tridegree at most $(1,1,1)$. We denote four of them by 
 $F_1,F_2,F_3,F_4$. 

 First, we assume that the $F_1$ is irreducible.
 The resultants of $F_1$ and $F_i$, 
 $i=2,3,4$ with respect to the last variable $t_3$ are denoted by $F_{12},F_{13},F_{14}$.
 The bidegrees of them are at most $(2,2)$.
 All these polynomials vanish on $K_{sym}$. If one of them is $0$, such as $F_{12}=0$, 
 then $F_1$ and $F_2$ must have a non-trivial common factor.
 This can only be $F_1$, since $F_1$ is irreducible. Then the tridegree of $F_1$ is less
 then $(1,1,1)$. Because $F_1$ vanishes on the non-degenerate component $K_0$, it must contain at least
 two variables, and so $F_1$ is a polynomial of degree $(1,1,0)$, maybe after some
 permutation of variables.

 If none of the three resultants vanishes, then let $G=gcd(F_{12}, F_{13}, F_{14})$. 
 The bidegree of $G$ is in the set $\{(2,2),\ (2,1),\ (1,1)\}$, up to permutation
 of variables $t_1,t_2$. 
 If it is $(1,1)$, then $G$ can be considered as a polynomial of tridegree $(1,1,0)$
 that vanishes on $K_0$.
 If the bidegree of $G$ is $(2,2)$ or $(2,1)$, then we write
 $F_{12}=G U_2,F_{13}=G U_3,F_{14}=G U_4$ with suitable polynomials $U_2,U_3,U_4$.
 The bidegrees of $U_2,U_3,U_4$ are at most $(0,1)$, hence
 $U_2,\ U_3, \ U_4$ are linear dependent, which means that there are three real number
 $\lambda_2,\ \lambda_3,\ \lambda_4$ such that 
 $$\lambda_2 F_{12}+ \lambda_3 F_{13}+ \lambda_4 F_{14}=0.$$
 As a consequence, we have
 $$Res(F_1, \lambda_2 F_{2}+ \lambda_3 F_{3}+ \lambda_4 F_{4}) =0,$$
 where $Res$ denotes the resultant.
 Then we can continue as in the case $F_{12}=0$ above. Again we get a polynomial of degree $(1,1,0)$, maybe after some
 permutation of variables.

Second, if $F_1$ is reducible, then it has two factors with degree $(1,1,0)$ and $(0,0,1)$,
 up to permutation of variables $t_1,t_2,t_3$. Again, $F_1$ vanishes on the non-degenerate component $K_0$, and so
 it must contain at least
 two variables, and so it is a polynomial of degree $(1,1,0)$, maybe after some
 permutation of variables.

In all cases above, we have a polynomial of tridegree $(1,1,0)$ vanishing on $K_0$.
 Since $\infty^3$ is in $K_{sym}$, it is of the form $bt_1+ct_2+d=0$, with $b,c,d\in\R$ and $bc \neq 0$,
 as stated in the lemma. We can use it to eliminate $t_2$: on $K_0$, we have $t_2=-\frac{b t_1+d}{c}$.

 The equations for the projection of $K_0$ to the $(t_1,t_3)$-plane
 can be obtained by substituting. We get the equation $\bf{ N^{\dagger} } \bf{X'}=0$, where
 ${\bf N^{\dagger}}:={\bf M^{\dagger} }{\bf L}$, and
\begin{equation*}\label{F22}
{\bf L}=\left[
 \begin{array}{cccccc}
  \frac{-b}{c}&            0& \frac{-d}{c}&            0&            0\\
             0&            1&            0&            0&            0\\
             0& \frac{-b}{c}&            0& \frac{-d}{c}&            0\\
             0&            0&            0&            1&            0\\
             0&            0& \frac{-b}{c}&            0& \frac{-d}{c}\\
             0&            0&            1&            0&            0\\
             0&            0&            0&            0&            1
\end{array}
 \right].             
\end{equation*}
 This follows from the fact that on $K_0$, we can replace $\bf{X}$ by ${\bf L}{\bf X'}$.
 Because $\operatorname{rank}(L)=5$, we also get
 $\operatorname{rank}({\bf N^{\dagger}})\geq \operatorname{rank}({\bf M^{\dagger}})-2$.
\qed
\end{pf}

\begin{pot}
  $r\geq 2$: 
 Assume, indirectly, that $r\le 1$.
 Then the system $ \bf{ M^{\dagger} } X=0$ is equivalent to zero or only one single equation in three variables, and 
 it will have at least a two-dimensional complex configuration set, which contradicts our assumption. 

 $r \leq 4$: Assume, indirectly, that $r\geq 5$. 
 Then from Lemma \ref{lm1}, the projection of 
  $K_{sym}$ to $(t_1,t_3)$ is defined by 
\begin{equation}\label{eq01}
 \bf{ N^{\dagger} }{\bf X}'=0,
\end{equation}
 where $r_{1}:=\operatorname{rank}({\bf{ N^{\dagger} }})\geq r-2 \geq 3$. 
 The equation (\ref{eq01}) is equivalent to a system of $r_{1}$ polynomial equations of bidegree at most $(2,1)$.  
 Because $K_{sym}$ is a curve and has non-degenerate components, the $r_{1}$ polynomials have a common factor with bidegree 
 at least $(1,1)$.
 Then $r_{1} \leq 2$ which contradicts
   to $r_{1}\geq 3$. 
\end{pot}

\section{Classification}

This section contains three parts. First, we show that the existence of a line symmetry implies $r=2$.
 Second, we show that $r=2$ or $r=3$ implies a line symmetry or another geometric consequence
 which we call the ``parallel property''. Third, we relate the case $r=4$ to a linkage described
 in \cite{HSScubep, part11}.

\subsection{Line Symmetric Linkages}
 
We now describe line symmetric 6R linkages in terms of dual quaternions.
 A 6R linkage $L=[h_1,h_2,h_3,h_4,h_5,h_6]$ is  line symmetric if and only if there 
 is a line represented by a dual quaternion $l$ such that $l^2=-1$ and
\begin{equation}\label{ls1}
 h_1=lh_4l^{-1},\ \quad h_2=lh_5l^{-1}, \ \quad h_3=lh_6l^{-1},
\end{equation}
 where $ll^{-1}=1$.
 Geometrically, the rotation around $l$ by the angle $\pi$ takes $h_i$ to $h_{i+3}$ for $i=1,2,3$.
\begin{proposition}\label{pro1}
 If $L$ is line symmetric, then $r=2$.
\end{proposition}
\begin{pf}
 As the norm of $l$ is equal to $1$, it follows $l^{-1}=-l$ 
 and we write (\ref{ls1}) as 
\begin{equation}\label{ls2}
 h_1=-lh_4l,\ \quad h_2=-lh_5l, \ \quad h_3=-lh_6l.
\end{equation}
 We define a map $\alpha$ from dual quaternion to itself as
 $$\alpha: \mathbb{DH} \longrightarrow \mathbb{DH},
 \quad\quad\quad\ h \ \longmapsto h+l\bar{h}l,$$
 where $\bar{h}$ denotes the conjugate of $h$ in dual quaternion. 
 It is true that all entries 
 of $M^{\dagger}_1$ are in $\operatorname{Im}(\alpha)$. For instance, we have
 $\alpha(h_1)=h_1-lh_1l=h_1+h_4=g_1,
 \alpha(h_5h_4)=h_5h_4+lh_4h_5l=h_5h_4-(lh_4l)(-lh_5l)=h_5h_4-h_1h_2,
 \alpha(h_6h_5h_4)=h_6h_5h_4-lh_4h_5h_6l=h_6h_5h_4+(-lh_4l)(-lh_5l)(-lh_6l)=h_6h_5h_4+h_1h_2h_3.$
 It is not difficult to prove that $\alpha$ is a $\mathbb{R}$-linear map. 
 If we  consider 
 ${ M^{\dagger}_1 }$ to be an $8\times7$ matrix with real entries,
 then $r_{2}:=\operatorname{rank}(M^{\dagger}_1)$ is less or equal to the dimension of $\operatorname{Im}(\alpha)$.
 W.l.o.g. we assume $l= \qi$. 
 We compute $\mathrm{Im}(\alpha)$ as 
 $\alpha({1})=1+\qi\qi=1-1=0, \alpha({\eps})=\eps+\eps \qi\qi=0,
  \alpha({\qi})={\qi}-\qi\qi\qi=2\qi, 
  \alpha(\qj)=\qj-\qi\qj\qi=0, 
  \alpha(\qk)=\qk-\qi\qk\qi=0,
  \alpha({\eps}{\qi})={\eps}{\qi}-{\eps}\qi\qi\qi=2{\eps}\qi, 
  \alpha({\eps}\qj)={\eps}\qj-{\eps}\qi\qj\qi=0, 
  \alpha({\eps}\qk)={\eps}\qk-{\eps}\qi\qk\qi=0$.
 Therefore, the dimension of $\operatorname{Im}(\alpha)$ is $2$.
 So we have $r_{2}\leq2$.
 
 The next step is to prove that all ${ M^{\dagger}_i }$ for $i=1,...,6$ are equal.
 It is true that the first three columns are equal in all ${ M^{\dagger}_i}$ for $i=1,...,6$.
 As $\operatorname{Im}(\alpha)$ is equal to  $\langle \qi, \eps \qi \rangle_{\mathbb{R}}$  and
 $g_1,g_2,g_3, h_6h_5-h_2h_3\in \operatorname{Im}(\alpha)$, we obtain
\begin{equation}\label{ls5}
 g_1\times g_2=g_1\times g_3=g_2\times g_3=(h_6h_5-h_2h_3)\times g_1=0,
\end{equation}
where $g\times h$ denotes the cross product
 of purely vectorial dual quaternions $g,h$.
 The equalities $M^{\dagger}_1=\dots=M^{\dagger}_6$ can be shown from (\ref{ls5}). For instance,
 $h_5h_4-h_1h_2-(h_1h_5-h_2h_4)=h_5\times h_4-h_1\times h_2-h_1\times h_5+h_2\times h_4=
 g_2\times h_4-h_1\times g_2=g_2\times g_1=0, h_1h_5-h_2h_4-(h_4h_2-h_5h_1)
 =h_1h_5-h_2h_4+(\overline{h_1h_5-h_2h_4})=0$ or
 $h_6h_5h_4+h_1h_2h_3 -( h_1h_6h_5+h_2h_3h_4)=-{\langle h_6, h_5 \rangle }h_4+{\langle h_2, h_3 \rangle }h_4
 -{\langle h_2, h_3 \rangle }h_1+{\langle h_6, h_5 \rangle }h_1 +(h_6\times h_5)\times h_4+h_1\times (h_2\times h_3)
 -h_1\times (h_6\times h_5) -(h_2\times h_3)\times h_4=(h_6\times h_5+h_3\times h_2)\times g_1=(h_6h_5-h_2h_3)\times g_1=0$,
 where $\langle g, h \rangle$ denotes the inner product of purely vectorial dual quaternions $g,h$.
 As a consequence, we have $r=r_{2}\leq2$.
 But we have $r\geq 2$ by Theorem \ref{Th1}, so $r=2$.
\qed
\end{pf}

\begin{remark}
 The well-known fact that line symmetric linkages are movable can also be obtained as
 a corollary from Theorem~\ref{pro1}. When $r=2$, then the configuration set is defined
 by 2 equations in 3 variables.
\end{remark}

\subsection{Linkages with Rank $2$ and $3$}

 In this subsection, we show that $r=2\ \mathrm{or}\ 3$ implies either a line symmetry
 or another property, defined as follows.
 We say that $L=[h_1,\dots,h_6]$ has the parallel property if $h_1 \parallel h_4$, $h_2 \parallel h_3$, $h_5 \parallel h_6$,
 maybe after some cyclic permutation of indices. In this section, we always
 assume that the rank of the $\lambda$-matrix of $L$ is $2$ or $3$.

 In the following, we use the technique of generic points of algebraic curves. 
 This simplifies the analysis a lot. Let $C$ be an irreducible algebraic curve. 
 Let $F$ be a field such $C$ can be defined by equations over $F$ (for instance $F=\Q$).
 Following \cite[Section 93]{W50}, we say that some point $p\in C$ is generic if it fulfills no 
 algebraic conditions defined by polynomials with coefficients in $F$, excerpt those that are a
 consequence of the equations of $C$. The existence of generic points is shown in \cite[Section 93]{W50};
 typically, the coordinates of a generic point are transcendental numbers.

 Let $K_0 \subset K_{sym}^+$ be an irreducible non-degenerate component of the linkage
 $L=[h_1,\dots,h_6]$, and let $\tau_0=(t'_1,t'_2,t'_3)$ be a generic point of $K_0$. 
 The configuration  $\tau_0$ corresponds to a set of rotations around the joint axes.
 When we apply these rotations, we get new positions for the $6$ lines, and we define the
 transformed linkage by $L'=[h'_1,h'_2,h'_3,h'_4,h'_5,h'_6]$. Note that $L$ and $L'$
 represent really the same linkage, just in different initial positions.

\begin{lemma}\label{lm2}
 If $\primal(g'_1)= 0$, then $L$ has the parallel property. Here $\primal(h)$ denotes the primal part of the dual quaternion $h$.
  More precisely, we will have 
 $h_1 \parallel h_4$, $h_2 \parallel h_3$, $h_5 \parallel h_6$, in all configurations in $K_0$.
\end{lemma}
\begin{pf}
 Assume that $\primal(g'_1)=0$.
 The parallelity of the first and fourth axis can be expressed as a set of polynomial 
 equations in the configuration parameters $(t_1,t_2,t_3)$. These equations are 
 fulfilled for the generic point $\tau_0$. By a well-known property of generic points it follows that they
 are fulfilled for all points in $K_{0}$. For this reason, the first and fourth axis are parallel at all position.
 
 Let $S=[p_1,p_2,p_3,p_4,p_5,p_6]$, where $p_i=\primal(h'_i)$ for $i=1,\ldots,6$.
 Then $S$ is a spherical linkage with the first and fourth axis coinciding at all positions. We can separate $S$ into two 3R linkages 
 $S_1=[p_1,p_2,p_3]$ and $S_2=[p_4,p_5,p_6]$.
 A 3R linkage is necessarily degenerate: either some angles are constant or some axes coincide.
 Since $t_2$ is not a constant in $K_0$, we obtain $p_2=\pm p_3$ or $p_1=\pm  p_2$.
 Since $t_3$ is not a constant in $K_0$, we obtain $p_2= \pm  p_3$ or $p_1= \pm  p_3$.
 If $p_2\ne \pm p_3$, then we have $p_1 =\pm p_2$ and $p_1= \pm  p_3$, a contradiction.
 So we obtain $p_2= \pm  p_3$.
 Similarly, we also have $p_5= \pm p_6$.
 
 Therefore, we get a linkage with $h'_1 \parallel h'_4$, $h'_2 \parallel h'_3$, $h'_5 \parallel h'_6$. 
 Since the parallel property is fulfilled for the generic point of 
 the configuration curve, it is fulfilled for all 
 points in $K_0$. In particular, the original linkage $L$ has the parallel property. 
\qed
\end{pf}

There is no $i$ such that $g'_i=0$ for $i=1,2,3$, because if $g'_i=0$ would be true, then the lines $h'_i$ and $h'_{i+3}$
 would be equal; the initial configuration was chosen generically, so the lines $h_i$ and $h_{i+3}$ would be equal 
 for all configurations in $K_0$, and this is not possible.
 Moreover, 
 it is not possible that two of $g_i$ for $i=1,2,3$ have 0 primal parts.
 In order to prove this, 
 we assume indirectly $\primal(g'_2)=0$ and $\primal(g'_3)=0$. 
 By Lemma \ref{lm2}, we get 
 $h_2 \parallel h_5$, $h_3 \parallel h_4$, $h_1 \parallel h_6$ and $h_3 \parallel h_6$, $h_4 \parallel h_5$, $h_1 \parallel h_2$.
 It follows that $L$ is a planar 6R Linkage which has mobility more than one.

Before the main theorem, we give several lemmas in the following.

\begin{lemma}\label{lm3}
  Let $a,b$ be two purely vectorial dual quaternions. If $a\times b=0$, then there is a dual number $\alpha$ such that
  $b=\alpha a$ or $a=\alpha b$, or the primal parts of $a$ and $b$ both vanish.
\end{lemma}
\begin{pf}
 Straightforward.
\qed
\end{pf}

In the next two proofs, we use the following argument from linear algebra. 
 Let $1\le i_1<\dots<i_r<i_{r+1}<\dots<i_s\le 7$ be integers.
 Let $A:=a_1M^{\dagger}_1+\dots+a_6M^{\dagger}_6$ be some
 linear combination of the matrices $M^{\dagger}_1,\dots,M^{\dagger}_6$, where $a_1,\ldots,a_6\in\R$.
 If the vector space generated by the columns $(i_1,\dots,i_s)$ of $M^{\dagger}$
 is already generated by the columns $(i_1,\dots,i_r)$ of $M^{\dagger}$, then the 
 vector space generated by the columns $(i_1,\dots,i_s)$ of $A$
 is also generated by the columns $(i_1,\dots,i_r)$ of $A$.

\begin{lemma}\label{lm4}
  If $g'_3\times g'_1=g'_2\times g'_1=0$, then $g'_2\times g'_3=0$.
\end{lemma}
 \begin{pf}
We distinguish two cases.

Case~I: $\primal(g'_1)\ne 0$. By Lemma~\ref{lm3}, there exist $\alpha_2,\alpha_3\in\D$ such that
 $g'_2=\alpha_2 g'_1$ and $g'_3=\alpha_3 g'_1$, and it follows that $g'_2\times g'_3=0$.

Case~II: $\primal(g'_1)=0$. Then $\primal(g'_2)\ne 0$ and $\primal(g'_3)\ne 0$. 
 If there exists $\alpha\in\D$ such that $g'_3=\alpha g'_2$, then $g'_2\times g'_3=0$.
 Otherwise, $g'_1$ is a dual multiple of $g_2'$ but $g'_3$ is not, so $g'_1,g'_2,g'_3$ are linearly independent.
 Then the first three columns generate the column space of $M^\dagger$. By linear algebra,
 the first three columns of $A:=M^{\dagger}_1+M^{\dagger}_4-M^{\dagger}_3-M^{\dagger}_6$ also
 generate the column space of $A$. But
\begin{equation}\label{rt01}
 A=[0,0,0,0, 2g'_3 \times g'_1, 2 g'_3\times g'_2,\ast]
\end{equation}
(we do not care about the last entry denoted by $\ast$), and it follows that $g'_2\times g'_3=0$.
\qed
 \end{pf}

\begin{lemma}\label{lm5}
 We have $g'_3\times g'_1=g'_2\times g'_1=g'_2\times g'_3=0$.
\end{lemma}
\begin{pf}
 Let $r_3$ be the dimension of the vector space generated by $g'_1,g'_2,g'_3$. If $r_3=1$, then
 it follows that $g'_3\times g'_1=g'_2\times g'_1=g'_2\times g'_3=0$. If $r_3=2$ or $r_3=3$,
 then the vector space $V$ generated by the first 6 columns of $M^\dagger$ is already generated
 by the first three and one of the other three columns.

Assume, for instance, that $V$ is generated by columns $(1,2,3,6)$. By linear algebra,
 the corresponding columns also generate the space of the first six columns of
 \[ M^{\dagger}_1+M^{\dagger}_4-M^{\dagger}_2-M^{\dagger}_5=[0,0,0,2 g'_2\times g'_1, 2g'_3 \times g'_1, 0,*] . \]
 This implies $g'_3\times g'_1=g'_2\times g'_1=0$, and by Lemma~\ref{lm4}, we also get $g'_2\times g'_3=0$.

If $V$ is generated by columns $(1,2,3,4)$, then the above linear algebra argument shows 
 $g'_1 \times g'_3 = g'_2\times g'_3=0$. The equality $g'_2\times g'_1=0$ follows again from by Lemma~\ref{lm4},
 applied to the linkage $[h_3,h_4,h_5,h_6,h_1,h_2]$. The third case, when $V$ is generated by columns $(1,2,3,5)$,
 is also similar.
\qed
\end{pf}

\begin{lemma}\label{lm6}
 If $\primal(g'_i)\ne 0$ for $i=1,2,3$, then $L'$ is line symmetric. 
\end{lemma}
  \begin{pf}
By Lemma~\ref{lm3}, there exists a dual quaternion $u$ and invertible dual numbers $\alpha_1,\alpha_2,\alpha_3$ such that
 $g'_i=\alpha_i u$ for $i=1,2,3$. Let $\beta:=u\bar{u}\in\D$. Because the primal part of $u$ is nonzero,
 the primal part of $\beta$ is positive, and $\frac{1}{\sqrt{\beta}}$ is defined. We set $l':=\frac{1}{\sqrt{\beta}}u$.
 Then $l'^2=-1$ and
 $g'_i h'_i ={h'_i}^2+h'_{i+3}h'_i={h'_{i+3}}^2+h'_{i+3}h'_1=h'_{i+3}g'_i$, hence $h'_{i+3}=g'_i h'_i {g'_i}^{-1}=l' h'_i l'^{-1}$
 for $i=1,2,3$.
\qed
\end{pf}

\begin{theorem}\label{Th2}
 If $r=2\ \mathrm{or}\ 3$, then $L$ has a  line symmetry or the parallel property.
\end{theorem}
\begin{pf} Let $K_0 \subset K_{sym}^+$ be an irreducible 
 non-degenerate component and $\tau_0=(t_1,t_2,t_3,t_1,t_2,t_3)$ be a generic point of $K_0$. 
 We get $L'=[h'_1,h'_2,h'_3,$ $h'_4,h'_5,h'_6]$ by applying the rotations specified in $\tau$.
 By Lemmas \ref{lm4}, \ref{lm5}, and \ref{lm6}, we conclude that $L'$ has a line symmetry or the parallel property.
 If a line symmetric linkage moves in an angle symmetric way, then the transformed linkage is
 also angle symmetric. This implies that when $L'$ is line symmetric, then $L$ is also line symmetric.
 On the other hand, if $L'$ has the parallel property, then parallelity holds for all points in $K_0$,
 in particular $L$ has the parallel property.
\qed
\end{pf}

\begin{theorem}\label{Th3}
 If $r=2$, then $L$ is line symmetric.
\end{theorem}
\begin{pf}
 By Theorem~\ref{pro1} and Theorem \ref{Th2}, we  may assume
 that $L$ has
 and parallel property 
 and $r=2$. 
 Let $L'=[h'_1,h'_2,h'_3,h'_4,h'_5,h'_6]$ be the linkage transformed by a generic position. 
 We may assume $h'_1 \parallel h'_4$, $h'_2 \parallel h'_3$,  $h'_5 \parallel h'_6$.
 The primal part of $g'_1$ is $0$ and the primal parts $g'_2$ and $g'_3$ are not.
 We define $l'$ as $\frac{1}{\sqrt{g'_2\overline{g_2'}}}g_2'$. Then $l'^2=-1$. By Lemma~\ref{lm5}, we also get
 $h'_2=-l'h'_5l'$ and $h'_3=-l'h'_6l'$ (see also the proof of Lemma~\ref{lm6}).
 Moreover, $g'_1$ is a real multiple of $\eps l'$, and $g'_1h'_1=h'_4g'_1$.
 By the last equation, the primal part of $h'_1+l'h'_4l'$ is zero. The dual part
 of $h'_1+l'h'_4l'$ is equal to $u:=g'_1-h'_4+l'h'_4l'$. The vectorial part of
 $ul'=g'_1l'-h'_4l'-lh'_4$ vanishes, so $u$ is a multiple of $l'$. On the other hand,
 the scalar product of $u$ with $l'$ also vanishes, hence $u=0$ and $h'_1=-l'h'_4l'$.
 It follows that $L'$ and $L$ are same line symmetric.
\qed
\end{pf}

In the end of this subsection, we give a construction of angle-symmetric 6R linkage with parallel property.
 The construction is based on the fact that we have a partially line symmetry taking $h_2$
 to $h_5$ and $h_3$ to $h_6$ (see Lemma \ref{lm3} and Lemma \ref{lm5} above).

\begin{const}\label{c1}{(Angle-Symmetric 6R Linkage with Parallel Property)}

 I.   Choose a rotation axis $u$ such that $u^2=-1$.  

 II.  Choose another rotation axis $h_1$ such that $h_1^2=-1$ and it is perpendicular to $u$.  

 III. Choose two parallel rotation axes $h_2$ and $h_3$  which are not perpendicular to $u$ such that $h_2^2=h_3^2=-1$. 

 IV.  Set $h_4=-uh_1u+r\eps u$, where $r$ is a random real number. 

 V.   Set $h_5=-uh_2u$ and $h_6=-uh_3u$. 

 VI.  Our angle-symmetric 6R Linkage with parallel property is $L=[h_1,h_2,$ $h_3,h_4,h_5,h_6]$. 
\qed
\end{const}

\begin{figure}
        \centering
        \begin{subfigure}[b]{0.3\textwidth}
                \centering
                \includegraphics[width=\textwidth]{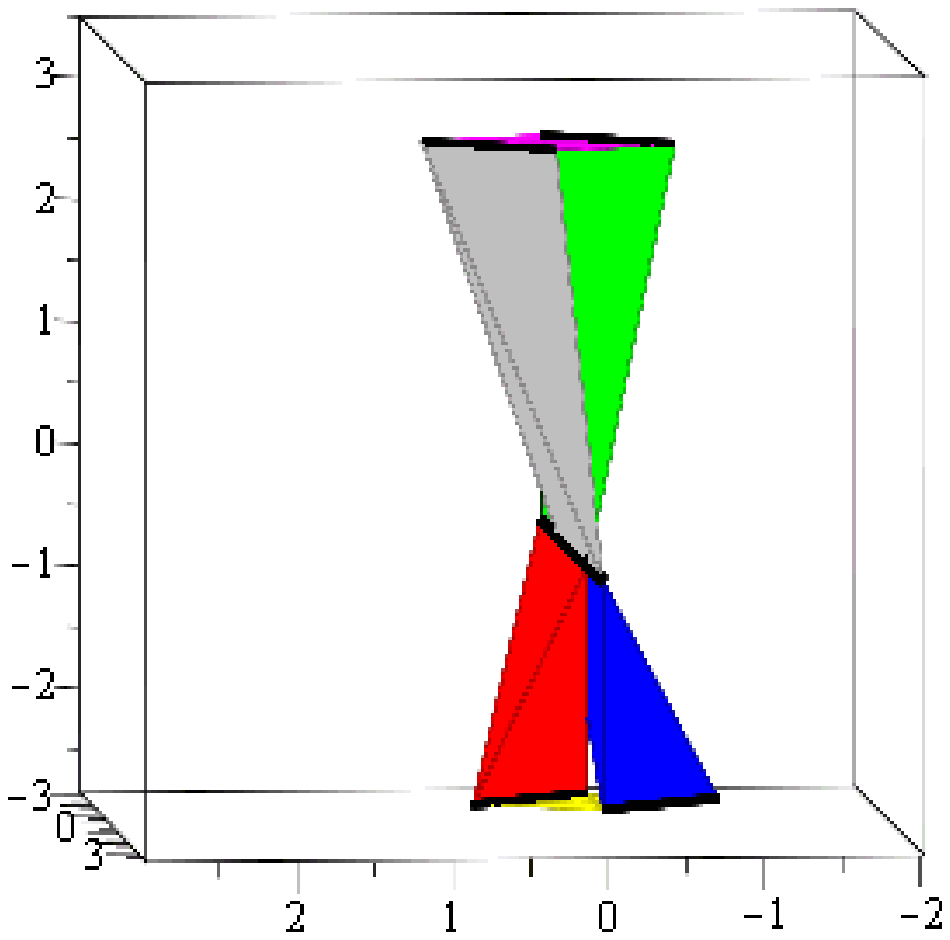}
                \caption{}
                \label{fig:1}
        \end{subfigure}
        \begin{subfigure}[b]{0.3\textwidth}
                \centering
                \includegraphics[width=\textwidth]{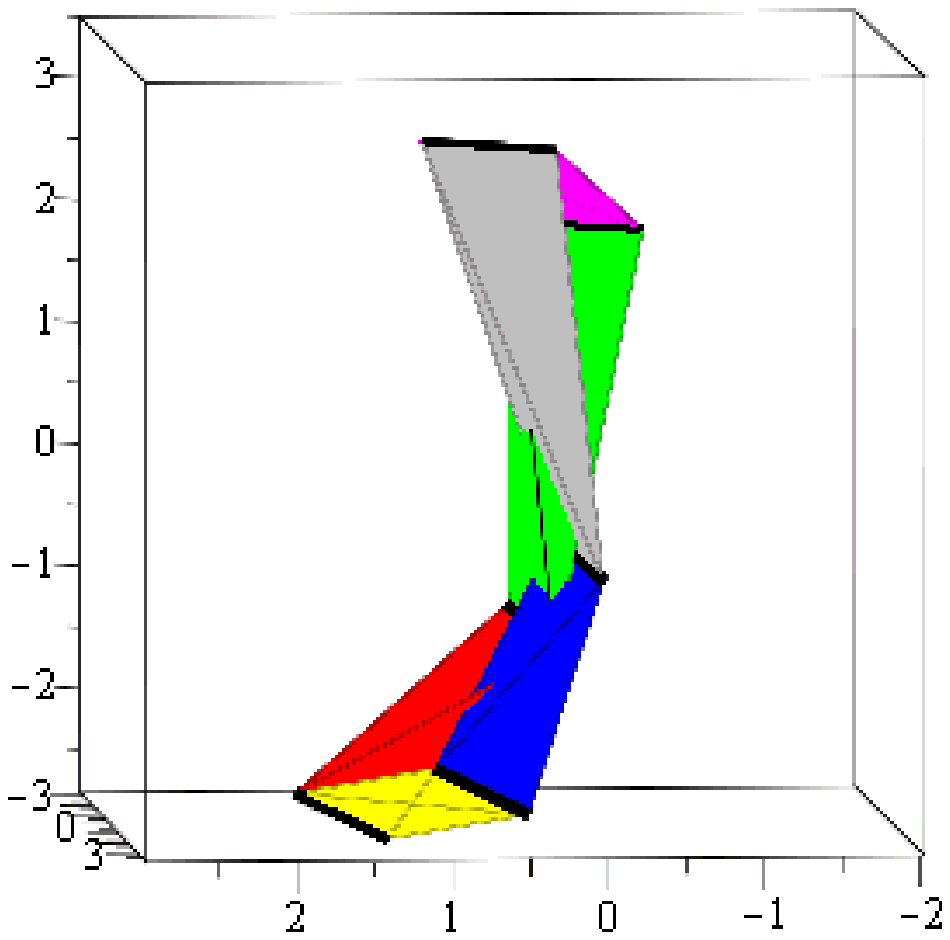}
                \caption{}
                \label{fig:2}
        \end{subfigure}
        \begin{subfigure}[b]{0.3\textwidth}
                \centering
                \includegraphics[width=\textwidth]{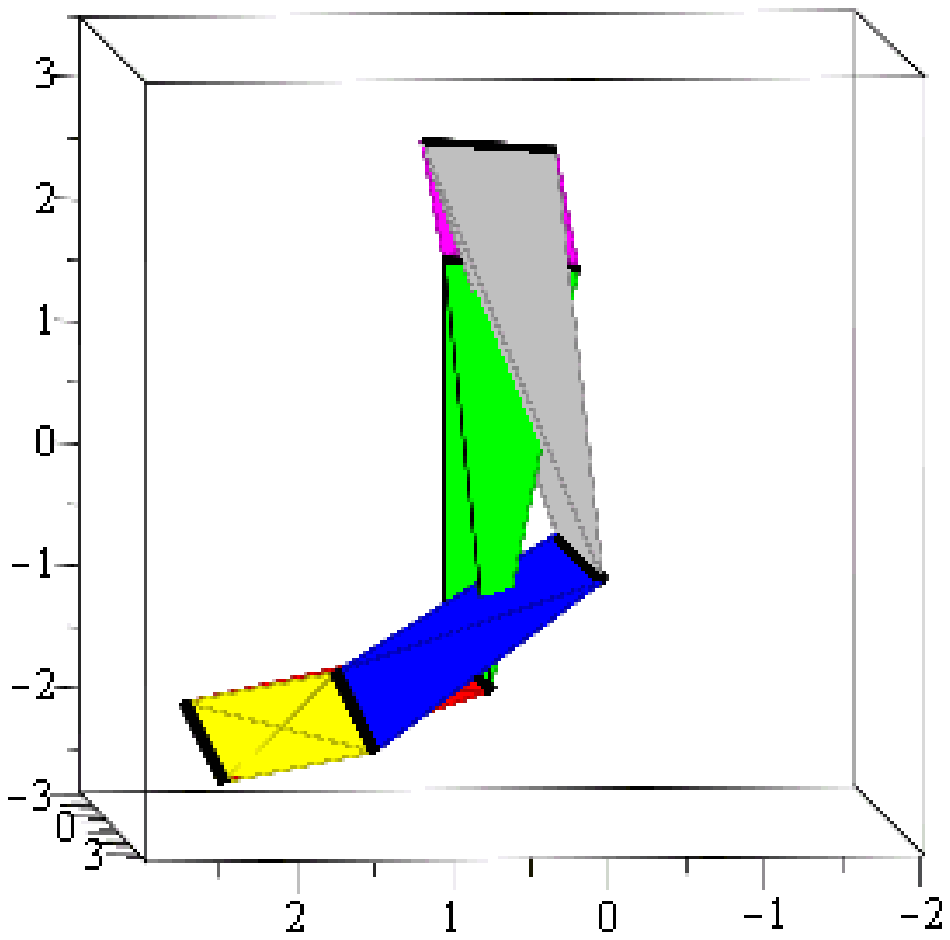}
                \caption{}
                \label{fig:3}
        \end{subfigure}

        \begin{subfigure}[b]{0.3\textwidth}
                \centering
                \includegraphics[width=\textwidth]{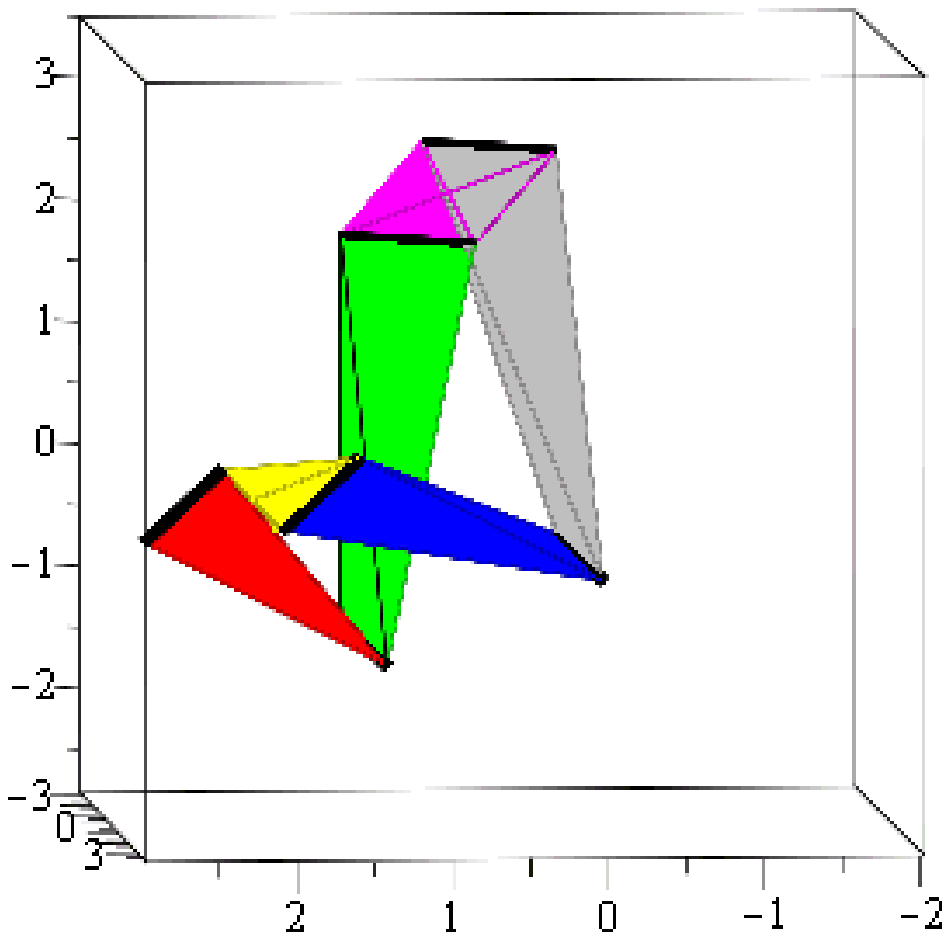}
                \caption{}
                \label{fig:4}
        \end{subfigure}
        \begin{subfigure}[b]{0.3\textwidth}
                \centering
                \includegraphics[width=\textwidth]{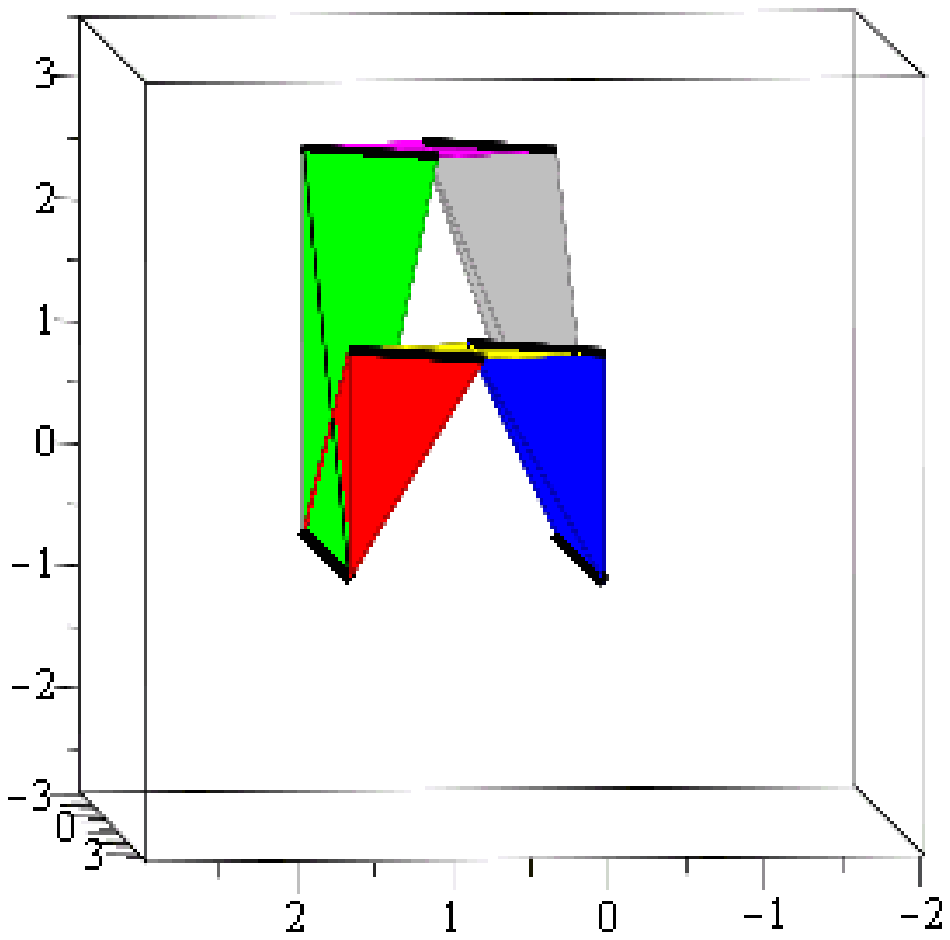}
                \caption{}
                \label{fig:5}
        \end{subfigure}
        \begin{subfigure}[b]{0.3\textwidth}
                \centering
                \includegraphics[width=\textwidth]{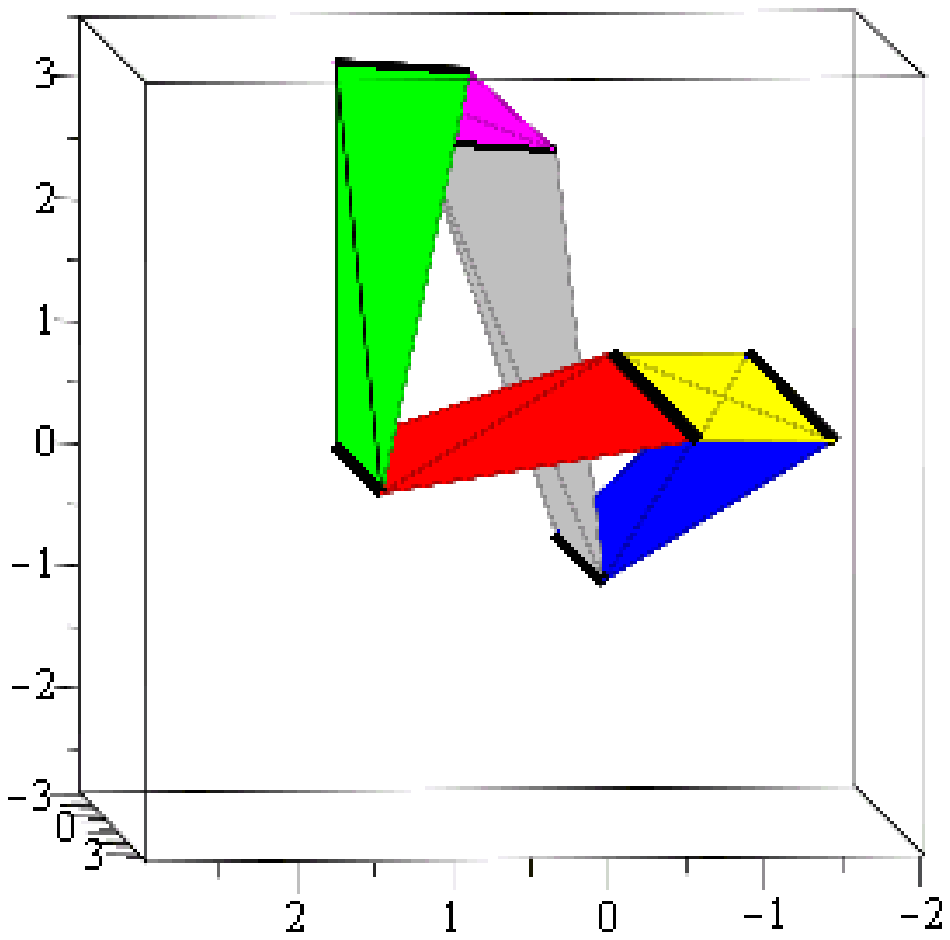}
                \caption{}
                \label{fig:6}
        \end{subfigure}
        \begin{subfigure}[b]{0.3\textwidth}
                \centering
                \includegraphics[width=\textwidth]{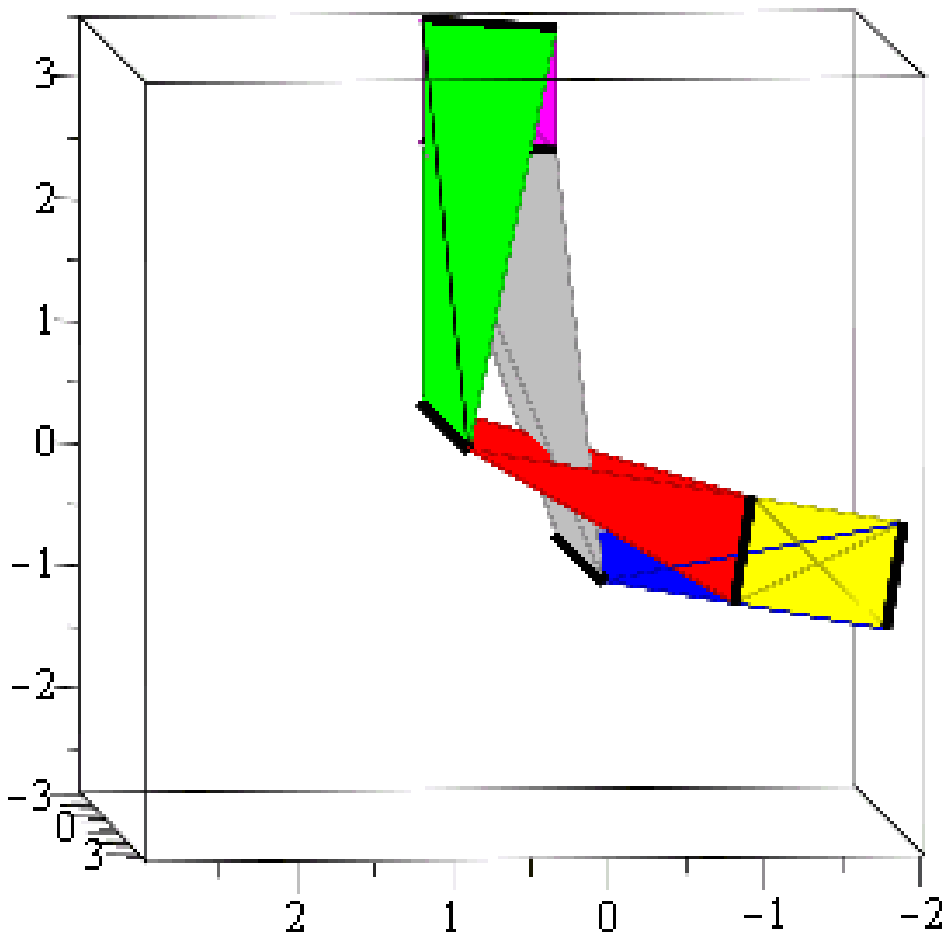}
                \caption{}
                \label{fig:7}
        \end{subfigure}
        \begin{subfigure}[b]{0.3\textwidth}
                \centering
                \includegraphics[width=\textwidth]{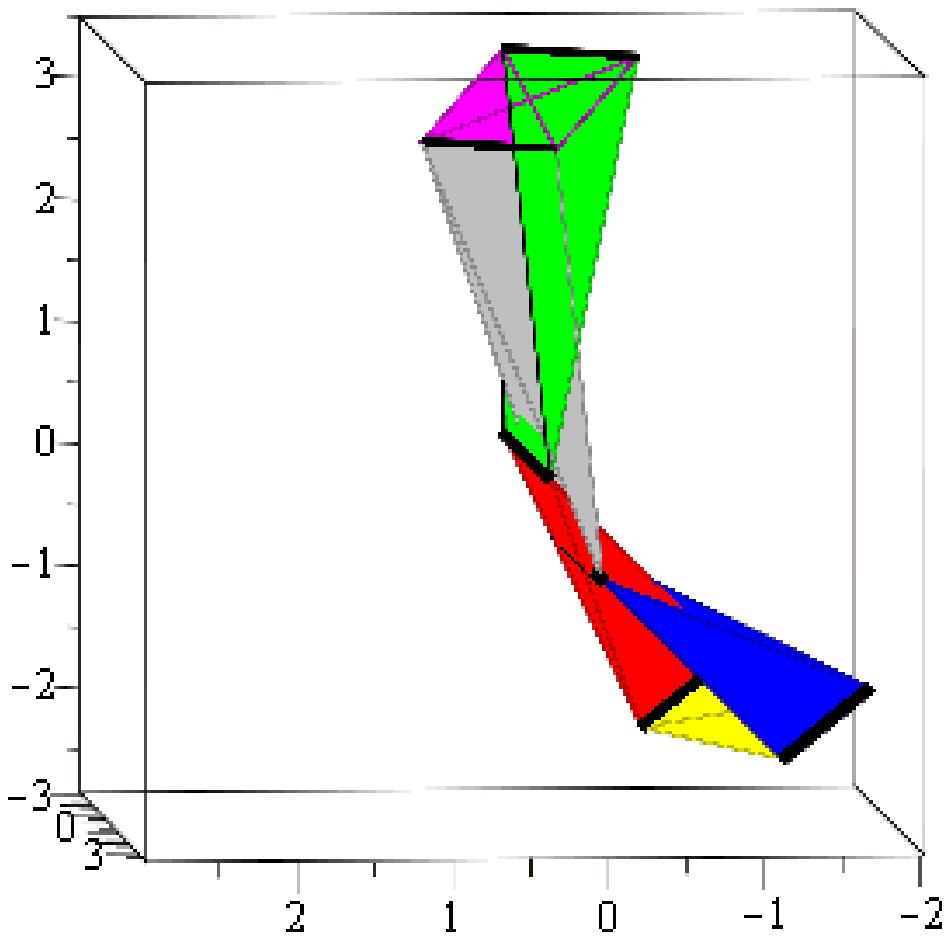}
                \caption{}
                \label{fig:8}
        \end{subfigure}
        \begin{subfigure}[b]{0.3\textwidth}
                \centering
                \includegraphics[width=\textwidth]{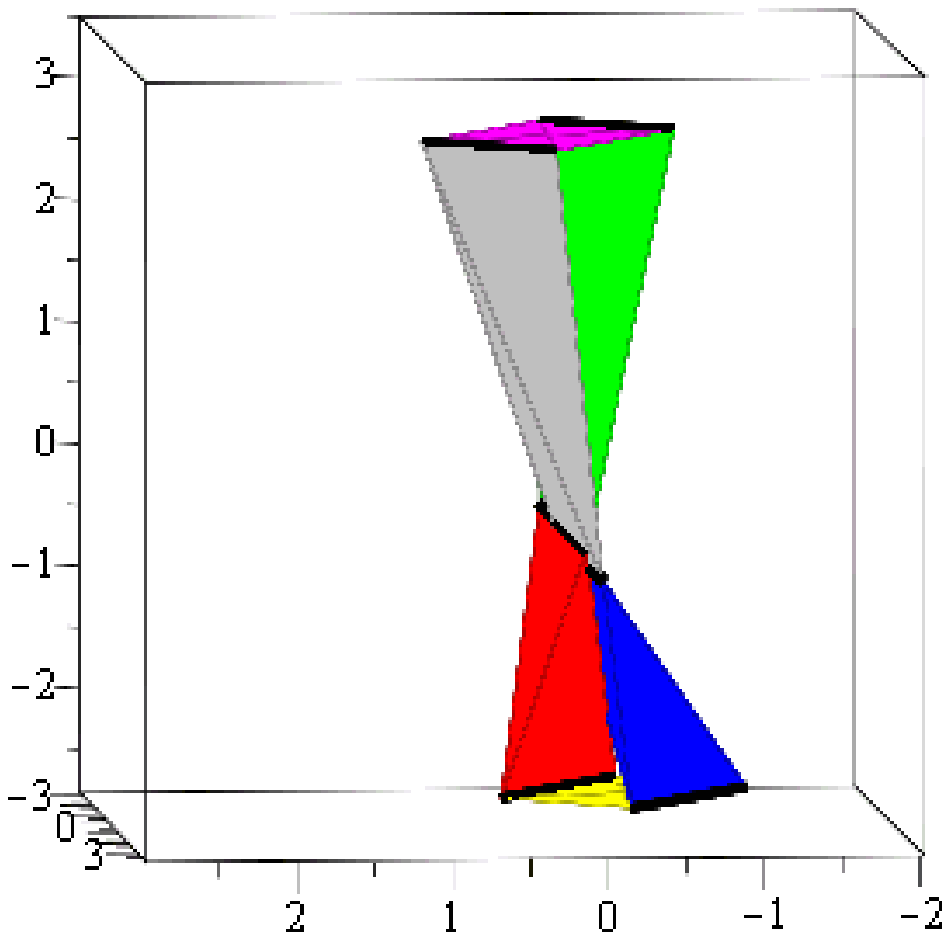}
                \caption{}
                \label{fig:9}
        \end{subfigure}
        \caption{These nine pictures which are produced by Maple are different positions of the linkage in Example \ref{ex1}.
 The six colored tetrahedra(gray, blue, yellow, red, green, pink) represent 
 six links in the linkage, and the joints are common edges of connected tetrahedra.}\label{fig1}
\end{figure}

\begin{example}{(Angle-Symmetric 6R Linkage with Parallel Property)}\label{ex1}
 We set 
 \begin{align*}
u&=\qi, \\
 h_{1}
&= -\frac{7}{11}\eps \qi +\qj, \\
h_{2}
&= \left(2\eps-\frac{3}{5}\right)\qi-
    \left(\frac{3}{2} \eps+\frac{4}{5}\right)\qj-
    \eps\qk, \\
h_{3}
&= \left(-2\eps+\frac{3}{5}\right)\qi+
    \left(\frac{3}{2} \eps+\frac{4}{5}\right)\qj+
    2\eps\qk, \\
 r&=\frac{14}{11},\\
 h_{4}
&= \frac{7}{11}\eps \qi -\qj, \\
h_{5}
&= \left(2\eps-\frac{3}{5}\right)\qi+
    \left(\frac{3}{2} \eps+\frac{4}{5}\right)\qj+
    \eps\qk, \\ 
 h_{6}
&= \left(-2\eps+\frac{3}{5}\right)\qi-
    \left(\frac{3}{2} \eps+\frac{4}{5}\right)\qj-
    2\eps\qk. \\
\end{align*}
It can be seen that the axes of $h_1$, $h_4$ are parallel, and the axes of $h_2,\ h_3$ and $h_5,\ h_6$,
 respectively, are parallel. Furthermore, the configuration curve contains a non-degenerate component:
\begin{align*}
(t_1,t_2,t_3,t_4,t_5,t_6)=\left(\frac{5}{4}t,t,t,\frac{5}{4}t,t,t\right).
\end{align*}
Thus, we have an example of angle-symmetric 6R linkage with parallel property. The rank of ${\bf M^{\dagger}}$ is 3.
 In Figure \ref{fig1}, 
 we present nine configuration positions of this linkage produced by Maple.
\qed
\end{example}

\begin{remark}
 A random instance of Construction \ref{c1} produces a linkage where $t_1$ is parametrized by a quadratic function in
 $t=t_2=t_3$. This example is special because $t_1$ is linear in $t$. (There is a degenerate component of the configuration curve
 that is responsible for this drop of the degree.)
\end{remark}

\subsection{Linkages with Rank $4$}

In this subsection, we show that the angle-symmetric linkages with Rank~4 are exactly those
that have been constructed in \cite[Example 3]{part11}  by factorization of cubic motion polynomials.

Recall that a motion polynomial $P$ is a polynomial in one variable $t$ with coefficients in $\Dh$
such that $P\bar{P}$ is a real polynomial that does not vanish identically. 
(Multiplication in $\Dh[t]$ is defined by requiring that $t$ commutes with the coefficients in $\Dh$.) 
Motion polynomials parametrize motions: by substituting a real number for $t$, we obtain an
element in the Study quadric. 

We give a brief sketch of the construction in \cite{HSScubep, part11}.
Linear motion polynomials of the form $(t-a-bh)$, $a,b\in\R$, $b\ne 0$, $h\in\Dh$, $h^2=-1$
parametrize revolutions. When we multiply three such polynomials $R_1,R_2,R_3$, we get a cubic
motion polynomial $Q$. Generically, there are 6 different factorizations into linear monic polynomials,
and there is one of the form $R_6R_5R_4$ such that the equations
$R_1\bar{R_1}=R_4\bar{R_4}$, $R_2\bar{R_2}=R_5\bar{R_5}$, $R_3\bar{R_3}=R_6\bar{R_6}$ hold. The three
linear factors $R_4,R_5,R_6$ are again motion polynomials parametrizing revolutions. The six axes
of $R_1,\dots,R_6$ define a closed 6R linkage; let us call it a linkage of cubic polynomial type.

We set $R_i(t)=t-a_i-b_ih_i$ for $i=1,\dots,6$, $a_i,b_i\in\R$, $b_i\ne 0$, $h_i\in\Dh$, $h_i^2=-1$.
The equations above are equivalent to $a_i=a_{i+3}$ and $b_i^2=b_{i+3}^2$ for $i=1,2,3$.
We may even assume $b_i=-b_{i+3}$; if not, we replace $h_{i+3}$ and $b_{i+3}$ by $-h_{i+3}$ and $-b_{i+3}$.
We multiply $R_1R_2R_3=R_6R_5R_4$ by $\bar{R_4}\bar{R_5}\bar{R_6}$ and get that
\[ (t-a_1-b_1h_1)(t-a_2-b_2h_2)(t-a_3-b_3h_3)(t-a_1-b_1h_4)(t-a_2-b_2h_5)(t-a_3-b_3h_6) \]
is a real polynomial. This shows that the configuration curve is parametrized by
\[ (t_1,t_2,t_3,t_4,t_5,t_6)=\left(\frac{t-a_1}{b_1},\frac{t-a_2}{b_2},\frac{t-a_3}{b_3},
	\frac{t-a_1}{b_1},\frac{t-a_2}{b_2},\frac{t-a_3}{b_3}\right) . \]
In particular, the linkage of cubic polynomial type is angle symmetric.

Here is a converse of the above statement.

\begin{theorem}\label{Th4}
If $L$ is an angle-symmetric linkage such that the $\lambda$-matrix has rank $r=4$, then
$L$ is of cubic polynomial type.
\end{theorem}
\begin{pf}
By Lemma \ref{lm1}, there exist a polynomial of the form $bt_1+ct_2+d$ that vanishes on $K_{sym}$,
$b,c,d\in\R$, $bc\ne 0$, 
and the projection of $K_{sym}$ to $(t_1,t_3)$ is in the common zero set of two linear
independent polynomials of bidegree $(2,1)$. The equation of the projection is therefore a common
factor of these two equations and must have bidegree smaller than $(2,1)$. Since $K_{sym}$ has a
non-degenerate component, the common factor cannot be constant in $t_1$ or $t_3$, hence it has bidegree
$(1,1)$. Because $(\infty,\infty)$ is contained in the projection, the common factor has the form
$b't_1+c't_2+d'$ for $b',c',d'\in\R$, $b'c'\ne 0$. This allows to parametrize $K_{sym}$ with
linear functions
\[ (t_1,t_2,t_3) = \left(\frac{t-a_1}{b_1},\frac{t-a_2}{b_2},\frac{t-a_3}{b_3}\right) \]
for $a_1,\dots,b_3\in\R$, $b_1b_2b_3\ne 0$. Now the linkage can be reconstructed from the
two factorizations of the cubic motion polynomial
\[ (t-a_1-b_1h_1)(t-a_2-b_2h_2)(t-a_3-b_3h_3)=(t-a_3+b_3h_6)(t-a_2+b_2h_5)(t-a_1+b_1h_4), \]
so it is of cubic polynomial type.
\qed
\end{pf}

\section{Conclusion}  
 In the analysis of the case $r=3$, we obtained a new type of linkages (with parallel property 
 $h_1 \parallel h_4$, $h_2 \parallel h_3$, $h_5 \parallel h_6$). It is not clear from the paper if every linkage with 
 parallel property is angle-symmetric. We know that this is not the case. A complete analysis of
 linkages with parallel property will be the topic of a future paper.

\section{Acknowledgements} 
 We would like to thank G\'abor Heged\"us and Hans-Peter Schr\"ocker
  for discussion and helpful remarks.
 The research was supported by the Austrian Science Fund (FWF):
W1214-N15, project DK9.

\end{document}